\documentclass[a4paper,fleqn]{cas-sc}

\usepackage[authoryear,longnamesfirst]{natbib}
\newtheorem{theorem}{Theorem}

 \newdefinition{rmk}{Remark}
 \newproof{pf}{Proof}
 \newproof{pot}{Proof of Theorem \ref{thm2}}
 
 \newtheorem{pro}[theorem]{Proposition}

\def\tsc#1{\csdef{#1}{\textsc{\lowercase{#1}}\xspace}}
\tsc{WGM}
\tsc{QE}
\tsc{EP}
\tsc{PMS}
\tsc{BEC}
\tsc{DE}

\begin{document}
\let\WriteBookmarks\relax
\def\floatpagepagefraction{1}
\def\textpagefraction{.001}
\shorttitle{Variation of Surfaces}
\shortauthors{N. Mosadegh}

\title [mode = title]{The Variational Geometry of Surfaces in the Conformally Flat Space}
\tnotemark[1]



\author[1,3]{N. Mosadegh}
\cormark[1]
\fnmark[1]{}
\ead[ur1]{n.mosadegh@azaruniv.ac.ir}


\address[1]{Depertment of Mathematics,
 Azarbaijan Shahid Madani University,
Tabriz 53751 71379, Iran}

\author[2,4]{E. Abedi}[]

\fnmark[2]
\ead{esabedi@azaruniv.ac.ir}
\ead[URL]{esabedi@azaruniv.ac.ir}


\address[2]{Depertment of Mathematics,
 Azarbaijan Shahid Madani University,
Tabriz 53751 71379, Iran}



\cortext[cor1]{Corresponding author}


\begin{abstract}
 In this paper, close surfaces are considered in 3-dimensional harmonic conformally flat space in point of the variation. It is shown that if the conformal vector field be tangent to surface and the sign of the mean curvature does not change then surface is minimal. Also, it is determined that the critical point of mean curvature functional is homeomorphic to the sphere. Furthermore, the Euler-Lagrange equations associated to the mean curvature and Willmore functionals are determined.
\\
\end{abstract}

\begin{keywords}
Variational Methods \sep Conformally Flat Spaces\sep Minimal Surfaces\\
MSC:\\
53A10, 53C42\\
\end{keywords}

\maketitle

\section{Introduction}
The research in the theory of surfaces in the Euclidean space are studied focus on the basic properties and the mean tools. The elementary theory of surfaces are presented for example see \citet{do}. Historically, as an interesting material surfaces are the local solutions of one of the oldest problem in geometry, called, the isoperimetric problem in geometry. Because they are compatible with mathematical models in outline in which the physical system asks a situation of least energy, so the surfaces are great interest for both mathematicians and physicists. Also, we refer to \citet{Ra} for the constant mean curvature surfaces with boundary, \citet{O,Ni} for minimal surfaces.

In recent years, several papers have presented the results which studied on surfaces with constant either extrinsic curvature or intrinsic curvature. For example the studies on surfaces with constant Gaussian curvature and with constant extrinsic curvature have been worked by \citet{MO,Es}, respectively.

Indeed, expanding studies on the surfaces renewed interest to deal with the close surfaces which are immersed in the 3-dimensional conformally flat space in which the ambient manifold is equipped with two Riemannian metrics $\widetilde{g}= e^{\sigma} <,>$, where $<,>$ denotes the Euclidean one. Also, if in this space the conformal map satisfies $\Delta_{\widetilde{g}} \sigma=\mathsf{tr}_{\widetilde{g}}\nabla^2 \sigma=0$ then it is called harmonic conformally flat manifold. Considering both the intrinsic and extrinsic geometry of a surface in point of the variational view is the aim of this paper, where they inherit the conformal metric of the ambient 3-dimensional harmonic conformally flat space. Indeed, we make use of some topological characterizations such as the Gauss-Bonnet theory and the Euler number to show the new results concerning surfaces in sense of the variation in our fundamental object of study. Briefly, we come to the corresponding variation to the curvature $K$ and the mean curvature $H$ of surfaces, where surfaces are tangent to the conformal vector field, then up to now the following results obtain
\begin{itemize}
  \item If the sign of mean curvature $H$ of a close surface does not change in the 3-dimensional harmonic conformally flat space, it is a minimal one.
  \item Under an isometric immersion a close complete oriented surface in the 3-dimensional harmonic conformally flat space, is homomorphic to the sphere. 
\item  The Euler- Lagrange equation associated to the critical points of the Willmore functional holds
 \begin{eqnarray*}
 \Delta H+ H\big(\|\omega^{\sharp}\|^2-\lambda_1\lambda_2 +2H^2-K \big)=0,
 \end{eqnarray*}
 where, $\lambda_i$ is the eigenvalue of the Weingarten operator and the vector field $\omega^{\sharp}$ is correlated to 1-form $\mathsf{\omega=d\sigma}$, metrically.
\end{itemize}
The paper is organized as follow. First of all, by means of the preparation for our investigation in this object, we need to have an overview of the geometry of the ambient 3-dimensional harmonic conformally flat space $\widetilde{M}^3$ and determine some results concerning isometrically immersed close surfaces there. Then, the principal results show in the rest of the sections 3 for the close surfaces in the 3-dimensional harmonic conformally flat space, where surfaces are tangent to the conformal vector field. Indeed, by applying the Gauss-Bonnet theory and integrating of the mean curvature in sense of variation most of the results obtain.
\section{Preliminaries}

 Let $M$ be a smooth Riemannian manifold which is equipped with the two conformal Riemannian metrics $g_1$ and $g_2$ in which $g_2=exp(\sigma)g_1$. Then, the Levi-Civita connections $^{2}\nabla$ and $^{1}\nabla$ of $g_2$ and $g_1$ are related
\begin{eqnarray}
^{2}\nabla_X Y = {^1}\nabla_{X}Y+\frac{1}{2}\{ \omega(X)Y +\omega(Y)X-g_2(X, Y)\omega^{\sharp}\},
\end{eqnarray}
where $X, Y\in \chi(M)$ and $\mathsf{\omega= d\sigma}$ for any 1-form $\omega$ in which the conformal vector field $\omega^{\sharp} = {^2\nabla} \sigma$. Next, we consider the following $(0, 2)$ tensor field
\begin{eqnarray}\label{am.1.21}
B(X, Y)= ({^1}\nabla_X \omega)(Y)-\frac{1}{2}\omega(X)\omega(Y),
\end{eqnarray}
such that, the curvature tensor fields $^2R$ and $^1R$ of $g_2$ and $g_1$ are connected
\begin{eqnarray}\label{am.1.1}
{^2}{R}(X, Y)Z= {^1}R(X, Y)Z &+& \frac{1}{2}\{B(X, Z)Y- B(Y, Z)X + g_2(X, Z)\ \ {^2}\nabla_Y \omega^{\sharp}- g_2(Y,Z)\ \ {^2}\nabla_X \omega^{\sharp}\}\nonumber\\
&+&\frac{1}{4}\big(g_2(Y, Z)\omega(X)- g_2(X, Z)\omega(Y)\big)\omega^{\sharp} ,
\end{eqnarray}
where $X$, $Y$ and $Z$ are tangent to $M$, and $(M,g_2)$ is known a conformally flat manifold if $g_1$ denotes the Euclidean metric. Also, $M$ is called harmonic conformally flat manifold provided that $\Delta \sigma=0$, where $\Delta=\mathsf{tr}_{g_2}\nabla^2$ is the Laplacian-Beltrami operator.

\begin{pro}\label{am.1.6}
Let $(M,g)$ be a surface in the 3-dimensional conformally flat space $(\widetilde{M}^3, \widetilde{g})$. Let $\{e_1, e_2\}$ be a locall orthonormal frame field on $M$, where $M$ is tangent to the $\omega^\sharp$. Then,
\begin{itemize}
  \item $\omega(e_1)^2=\omega(e_2)^2 $ and $\omega^{\sharp}=\|\omega^\sharp\|\frac{\sqrt{2}}{2}(e_1+e_2)$.
  \item $\widetilde{K}=\big(\frac{1}{4}\|\omega^{\sharp}\|^2 -\frac{1}{2}\mathsf{div} \omega^{\sharp}\big)$.
\end{itemize}
where, $\mathsf{div}\omega^{\sharp}=\mathsf{trace}\widetilde{\nabla} \omega^{\sharp}$ is associated to $M$.
\end{pro}
\begin{pf}
From (\ref{am.1.1}), the sectional curvature of plane, which is spanned by $\{e_1, e_2\}$, holds
\begin{eqnarray}\label{am.1.4}
\widetilde{g}(\widetilde{R}(e_1,e_2)e_2,e_1)=\frac{1}{2}\big(\|\omega^{\sharp}\|^2-\mathsf{div} \omega^{\sharp}-\omega(e_2)^2 \big),
\end{eqnarray}
analogously,
\begin{eqnarray}\label{am.1.5}
\widetilde{g}(\widetilde{R}(e_2,e_1)e_1,e_2)=\frac{1}{2}\big(\|\omega^{\sharp}\|^2 - \mathsf{div} \omega^{\sharp}-\omega(e_1)^2\big),
\end{eqnarray}
 Then, immediately these equations follow $\omega(e_1)^2=\omega(e_2)^2$. From here we can have $\theta=\frac{\pi}{4}$, where $\theta$ is angle between $\omega^{\sharp}$ and the vector field $e_i$. Then
\begin{eqnarray}
\omega^{\sharp}=\widetilde{g}(\omega^{\sharp},e_1)e_1+\widetilde{g}(\omega^{\sharp},e_2)e_2=\frac{\sqrt{2}}{2}\|\omega^{\sharp}\|(e_1+ e_2).
\end{eqnarray}
 Next, adding (\ref{am.1.4}) and (\ref{am.1.5}) yield
\begin{eqnarray}\label{am.1.8}
\widetilde{K}=\frac{1}{2}\big(\frac{1}{2}\|\omega^{\sharp}\|^2 -\mathsf{div} \omega^{\sharp}\big).
\end{eqnarray}
So, we reach the result.
\end{pf}
Now, we restrict our attention to a Riemannian close surface $M$ in the 3-dimensional conformally flat space $\widetilde{M}$, where $M$ is tangent to the $\omega^{\sharp}$. Consider $\widetilde{\nabla}$ and $\nabla$ are the Levi-Civta connections on $\widetilde{M}$ and $M$, as well. Let $e_i$ , $N$  and $H$ be the principal direction associated to the eigenvalues $\lambda_i$, $i=1,2$, the local unit normal vector field and the mean curvature on $M$, respectively. By $\widetilde{K}$ we denotes the sectional curvature spanned by $\{e_1,e_2\}$ on $\widetilde{M}$ and $K$ is the Gaussian curvature on $M$.

Further the Codazzi equation for a surface $M$ in the 3-dimensional conformally flat space $(\widetilde{M}^3, \widetilde{g})$ holds
\begin{eqnarray}
\widetilde{g}\big( (\widetilde{\nabla}_X A)Y-(\widetilde{\nabla}_Y A)X, Z\big)=\frac{1}{2}\big(\omega(AY)\widetilde{g}(X, Z)-\omega(AX)\widetilde{g}(Y,Z)\big),\label{am.1.7}
\end{eqnarray}
where, $A$ is the Weingarten operator for vector fields $X, Y$ and $Z$ tangent to $M$.

\begin{pro}\label{am.1.10}
Let $M$ be a close surface in the 3-dimensional conformally flat space $(\widetilde{M}^3, \widetilde{g})$. Then
\begin{eqnarray}
  \int_M \lambda_2\widetilde{g}(\nabla_{e_1}X, e_1)+\lambda_1\widetilde{g}(\nabla_{e_2}X,e_2)-\frac{1}{2}\omega(AX)d\Omega=0,
\end{eqnarray}
for $X\in \chi(M)$.
\end{pro}
\begin{pf}
Let $e_i$ for $i=1,2$ be the eigenvectors of the Weingarten operator $A$ corresponding to the eigenvalues $\lambda_i$. By taking into account (\ref{am.1.7}), from the following calculation we have
\begin{eqnarray}\label{am.1.26}
\lambda_2\widetilde{g}(\nabla_{e_1}X, e_1)+\lambda_1\widetilde{g}(\nabla_{e_2}X,e_2)&=& 2H\mathsf{div} X-\big(\lambda_2\widetilde{g}(\nabla_{e_2}X, e_2)+\lambda_1\widetilde{g}(\nabla_{e_1}X,e_1) \big)\nonumber\\
&=&2H\mathsf{div} X-\big(\widetilde{g}(\nabla_{e_2}X, Ae_2)+\widetilde{g}(\nabla_{e_1}X, Ae_1) \big)\nonumber\\
&=&2H\mathsf{div} X-\big(\widetilde{g}(\nabla_{e_2}AX-(\nabla_{e_2}A)X, e_2)+\widetilde{g}(\nabla_{e_1}AX-(\nabla_{e_1}A)X,e_1) \big)\nonumber\\
&=&2H\mathsf{div} X- \mathsf{div} AX+\big(\widetilde{g}((\nabla_{e_2}A)X, e_2)+\widetilde{g}((\nabla_{e_1}A)X,e_1) \big)\nonumber\\
&=&2H\mathsf{div} X- \mathsf{div} AX +\big(\widetilde{g}((\nabla_X A)e_2+\frac{1}{2}(\omega(AX)e_2-\omega(Ae_2)X), e_2)\nonumber\\
&&+\widetilde{g}((\nabla_X A)e_1+\frac{1}{2}(\omega(AX)e_1-\omega(Ae_1)X),e_1) \big)\nonumber\\
&=&2H\mathsf{div} X- \mathsf{div} AX +\big(\widetilde{g}((\nabla_X A)e_2,e_2)+\widetilde{g}((\nabla_X A)e_1,e_1)\big)+\frac{1}{2}\omega(AX)\nonumber\\
&=&2H\mathsf{div} X- \mathsf{div} AX+ (X\lambda_1)+(X\lambda_2)+\frac{1}{2}\omega(AX)\nonumber\\
&=&2H\mathsf{div} X- \mathsf{div} AX+ 2\widetilde{g}(\nabla H, X)+\frac{1}{2}\omega(AX)\nonumber\\
&=&\mathsf{div}( 2HX)- \mathsf{div} AX + \frac{1}{2}\omega(AX)\nonumber,
\end{eqnarray}
for all $X\in \chi(M)$. Since, the surface is closed by applying the divergence theorem we conclude the proof.
\end{pf}

\begin{pro}\label{am.1.23}
Let $M$ be a close surface in the $3$-dimensional conformally flat space $(\widetilde{M}^3,\widetilde{g})$. Then
\begin{eqnarray}
  \int_M \lambda_2 f \mathsf{Hess}_g(e_1, e_1)+\lambda_1 f \mathsf{Hess}_g(e_2,e_2)-\frac{1}{2}f\omega (A \nabla g)d \Omega\nonumber\\
=\int_M \lambda_2 g \mathsf{Hess}_f(e_1, e_1)+\lambda_1 g \mathsf{Hess}_f(e_2,e_2)-\frac{1}{2}g\omega (A \nabla f)d \Omega,
\end{eqnarray}
where $f, g\in C^{\infty}(M)$.
\end{pro}
\begin{pf}
We are in the hypotheses of the last proposition. Let $X=(f\nabla g-g\nabla f) \in \chi(M)$, where $f, g \in C^{\infty}(M)$, then the divergence theorem states $\int_M f \Delta g d\Omega=\int_M g \Delta f d\Omega$. Next, Proposition \ref{am.1.10} follows
\begin{eqnarray}
0=\int_M \big\{\lambda_2\widetilde{g}(\widetilde{\nabla}_{e_1}(f\nabla g-g\nabla f),e_1)+\lambda_1\widetilde{g}(\widetilde{\nabla}_{e_2}(f\nabla g-g\nabla f),e_2)-\frac{1}{2}\omega(A(f\nabla g-g\nabla f))\big\}d\Omega,
\end{eqnarray}
then
\begin{eqnarray}
\int_M \big\{ f\lambda_2\widetilde{g}(\widetilde{\nabla}_{e_1}\nabla g,e_1)+f\lambda_1\widetilde{g}(\widetilde{\nabla}_{e_2}\nabla g,e_2)-\frac{1}{2}f\omega(A\nabla g)\big\}d\Omega,\nonumber\\
=\int_M \big\{ g\lambda_2 \widetilde{g}(\widetilde{\nabla}_{e_1}\nabla f,e_1)+g\lambda_1 \widetilde{g}(\widetilde{\nabla}_{e_2}\nabla f,e_2)-\frac{1}{2}g \omega(A\nabla f))\big\}d\Omega,
\end{eqnarray}
since $\widetilde{g}(\widetilde{\nabla}_{e_i}\nabla f, e_i)=\mathsf{Hess}_f(e_i,e_i)$ and analogously for $g$, then the proof is completed.
\end{pf}

 A variation of an isometric immersion $\varphi: M\rightarrow \widetilde{M}$ is considered a differentiable map $\phi : M \times (-\epsilon, \epsilon)\rightarrow \widetilde{M}: (p, t)\mapsto \varphi_t(p)=\phi(p,t)$, where $\varphi_t: M\rightarrow \widetilde{M}$ is an immersion such that for $t=0$, $\varphi_0=\varphi$. Furthermore, the variational vector field $\xi$, associated to $\{\varphi_t\}$ attaches $\xi_{p}=\frac{\partial \phi(p,t)}{\partial t}|_{t=0}$ for $p\in M$. The variation $\{\varphi_t\}$ is said to be a tangent variation where the vector field $\xi$ is tangent to the surface in any points. It is convenient to set convention the variation is normal where the variational vector field $\xi=f N$, where $N$ is the local unit normal vector field of $M$ and $f\in C^{\infty}(M)$.

In the rest of this section, we come up to establish the variation of the eigenvalues of the Weingarten operator inspired by the variation of the Weingarten operator that was given:
\begin{theorem}(\citet{V} Variation of the Weingarten operator)
Let $(M,g)$ be a semi-Riemannian hypersurface of a semi-Riemannian manifold $\widetilde{M}$ under a variation of $M$, with variational vector
field $\xi = f N +\xi^t$, then for every $X \in \chi(M)$ holds :
\begin{eqnarray}\label{am.1.24}
(\delta A)(X)=f\big(\widetilde{R}(N,X)N+A^2(X)\big)+\mathsf{Hs_f}(X)+(\mathsf{\pounds}_{\xi^t}A)(X),
\end{eqnarray}
where $\mathsf{Hs_f}(X)$ is the operator associated to the tensor of $\mathsf{Hess}_f(X,Y)=\widetilde{g}(\mathsf{Hs_f}(X), Y)= \widetilde{g}(\widetilde{\nabla}_X\nabla f, Y )$ and $(\pounds_{\xi^{t}}A)(X)$ is the Lei derivative of the Weingarten operator along $\xi^t$.
\end{theorem}
 Let $\varphi_t:M \rightarrow \widetilde{M}^3$ be an immersed normal variation of a close surface $M$ in the 3-dimensional conformally flat space $\widetilde{M}^3$ where $\xi=fN$. Let $e_i$ be the principal direction on $M$, corresponding to the eigenvalue $\lambda_i$ of the Weingarten operator $A$. Now, taken the variation of both sides of the equation $A(\varphi_t)e_i(\varphi_t)= \lambda_i(\varphi_t)e_i(\varphi_t)$, leads to
\begin{eqnarray}
(\delta A)(e_i)+A(\delta e_i)=(\delta \lambda_i)e_i+\lambda_i(\delta e_i),\nonumber
\end{eqnarray}
since $\widetilde{g}(\delta e_i,e_i)=0$ then $\delta e_i=\widetilde{g}(\delta e_i, e_j)e_j$ for $i,j=1,2$ and $j\neq i$, it follows
\begin{eqnarray}
(\delta A)(e_i)+ \lambda_j(\delta e_i)=(\delta \lambda_i)e_i+\lambda_i(\delta e_i), \ \ \ \ \ i\neq j\nonumber
\end{eqnarray}
then $\delta \lambda_i= \widetilde{g}((\delta A)e_i, e_i)$. From here, by applying $(\ref{am.1.24})$ we get
\begin{eqnarray}\label{am.1.2}
\delta \lambda_i=f\big(\widetilde{g}(\widetilde{R}(N, e_i)N, e_i)+\lambda^2_i\big)+ \mathsf{Hess}_f(e_i,e_i).
\end{eqnarray}

\section{Curvature Properties and Variational Aspects of surfaces in the 3-dimensional harmonic conformally flat space}
From now we concentrate on the main results.
\begin{theorem}\label{am.1.18}
Let $\varphi_t: M\rightarrow \widetilde{M}^3$ be an immersed normal variation of a close surface $M$ in the 3-dimensional harmonic conformally flat space. Then, the variation of the Guassian curvature of $M$ holds
\begin{eqnarray}
\delta K &=& \frac{f}{2}\big( (\lambda_1-\lambda_2) \widetilde{g}(\widetilde{\nabla}_{e_1}\omega^{\sharp}, e_1)+ (\lambda_2-\lambda_1)\widetilde{g}(\widetilde{\nabla}_{e_2} \omega^{\sharp},e_2)\big)\nonumber\\
&& + \lambda_1\mathsf{Hess}_f(e_2,e_2)+\lambda_2\mathsf{Hess}_f(e_1,e_1) -\frac{1}{2}\mathsf{div}(\delta \omega^{\sharp})\nonumber\\
 && + f\big( 2H(\lambda_1\lambda_2) - H\omega(\frac{\nabla f}{f})+\frac{1}{8}H\| \omega^{\sharp}\|^2 -\omega(\frac{A\nabla f}{f})\big).
\end{eqnarray}
\end{theorem}
\begin{pf}
Let $\varphi:M \rightarrow \widetilde{M}^3$ be an isometric immersion of a close surface $M$ into the 3-dimensional harmonic conformally flat space $\widetilde{M}$. Let $e_i$ and $N$ be the principal direction corresponding to the eigenvalues $\lambda_i$, $i=1,2$ and the local unit normal vector field on $M$, respectively. We consider the Guass equation that can be written as $K=\widetilde{K}+\lambda_1\lambda_2$. From here we have
\begin{eqnarray}\label{am.1.14}
\frac{\partial}{\partial t}\big|_{t=0}K=\delta K=\delta \widetilde{K}+\delta(\lambda_1\lambda_2).
\end{eqnarray}
We consider equations (\ref{am.1.21}) and (\ref{am.1.8}), impliy
\begin{eqnarray}\label{am.1.15}
\delta \widetilde{K}=-\frac{1}{2}\delta\big( B(e_1,e_1)+ B(e_2, e_2)\big).
\end{eqnarray}
such that
\begin{eqnarray}
\frac{\partial}{\partial t}\big|_{t=0}B(e_i,e_i)=\delta B(e_i,e_i)
= \delta\big( \widetilde{g}(\widetilde{\nabla}_{e_i}\omega^{\sharp}, e_i)-\frac{1}{4}\|\omega^{\sharp}\|^2 \big).\nonumber
\end{eqnarray}
The next step is to compute variation of the first and second terms of the last equation. We obtain
\begin{eqnarray}\label{am.1.31}
\delta\widetilde{g}(\widetilde{\nabla}_{e_i}\omega^{\sharp}, e_i)&=& (\delta \widetilde{g})(\widetilde{\nabla}_{e_i}\omega^{\sharp}, e_i)+ f\widetilde{g}(\widetilde{\nabla}_N\widetilde{\nabla}_{e_1}\omega^{\sharp}, e_i) + f\widetilde{g}(\widetilde{\nabla}_{e_i}\omega^{\sharp}, \widetilde{\nabla}_N e_i),\nonumber\\
&=&-2\lambda_1f\widetilde{g}(\widetilde{\nabla}_{e_i}\omega^{\sharp}, e_i)+ f \widetilde{g}(\widetilde{R}(N,e_i)\omega^{\sharp}, e_i) +\widetilde{g}(\widetilde{\nabla}_{e_i}\delta \omega^{\sharp},e_i)+ f\widetilde{g}(\widetilde{\nabla}_{[N, e_i]}, e_i)\nonumber\\
&&+f\mathsf{Hess}_{\sigma}(e_i,\widetilde{\nabla}_N e_i),
\end{eqnarray}
where $(\delta\widetilde{g})(X, Y)=-2f\widetilde{g}(AX, Y)$ for $X,Y\in \chi(M)$ see(\citet{V}). Also,
\begin{eqnarray}\label{am.1.30}
\delta \|\omega^{\sharp}\|^2 &=&(\delta\widetilde{g})(\omega^{\sharp},\omega^{\sharp})+\frac{1}{2}\widetilde{g}(\delta \omega^{\sharp},\omega^{\sharp}),\nonumber\\
&=&-2f\widetilde{g}(A\omega^{\sharp}, \omega^{\sharp})+\frac{1}{2}f\widetilde{g}(\widetilde{\nabla}_N\omega^{\sharp},\omega^{\sharp}),\nonumber\\
&=&-2f\omega(A \omega^{\sharp})+ \frac{f}{2}B(N, \omega^{\sharp})=-\frac{3}{2}f\omega(A\omega^{\sharp}),
\end{eqnarray}
where (\ref{am.1.21}) yields $B(N, \omega^{\sharp})=-\omega(\widetilde{\nabla}_{\omega^{\sharp}}N)=\omega(A \omega^{\sharp})$. Equations (\ref{am.1.15}) together with (\ref{am.1.31}) and (\ref{am.1.30}), lead to
\begin{eqnarray}\label{am.1.22}
\delta \widetilde{K}=
\frac{1}{2}f\big( \lambda_1 \widetilde{g}(\widetilde{\nabla}_{e_1}\omega^{\sharp}, e_1)+ \lambda_2\widetilde{g}(\widetilde{\nabla}_{e_2} \omega^{\sharp},e_2)\big)-\frac{5}{8}fH\| \omega^{\sharp}\|^2 -\frac{1}{2}\mathsf{div}(\delta \omega^{\sharp})-f\omega(\frac{A\nabla f}{f}),
\end{eqnarray}
because of
\begin{eqnarray}\label{am.1.17}
\omega(A\omega^{\sharp})&=&\widetilde{g}(\omega^{\sharp}, \widetilde{g}(\omega^{\sharp},e_1)Ae_1+\widetilde{g}(\omega^{\sharp},Ae_2)e_2)\nonumber\\
&=& \lambda_1\omega(e_1)^2+\lambda_2\omega(e_2)^2=H\|\omega^{\sharp}\|^2.
\end{eqnarray}
Further, equation $(\ref{am.1.2})$ notes that
\begin{eqnarray}\label{am.1.16}
\delta(\lambda_1\lambda_2)&=&  2fH(\lambda_1\lambda_2) + \lambda_1\mathsf{Hess}_f(e_2,e_2)+\lambda_2\mathsf{Hess}_f(e_1,e_1)\nonumber\\
&&+ f\lambda_2\big(-\frac{1}{2}(\widetilde{g}(\widetilde{\nabla}_{e_1}\omega^{\sharp}, e_1)+B(N,N))+\frac{1}{4}\omega(e_1)^2 \big)\nonumber\\
&&+f \lambda_1\big(-\frac{1}{2}(\widetilde{g}(\widetilde{\nabla}_{e_2}\omega^{\sharp}, e_2)+B(N,N))+\frac{1}{4}\omega(e_2)^2 \big),\nonumber\\
&=& \lambda_1\mathsf{Hess}_f(e_2,e_2)+\lambda_2\mathsf{Hess}_f(e_1,e_1)-\frac{f}{2}\big(\lambda_2\widetilde{g}(\widetilde{\nabla}_{e_1}\omega^{\sharp}, e_1)+\lambda_1\widetilde{g}(\widetilde{\nabla}_{e_2}\omega^{\sharp}, e_2)\big)\nonumber\\
 && + f\big(2H(\lambda_1\lambda_2) - H\omega(\frac{\nabla f}{f}) + \frac{3}{4} H\|\omega^{\sharp}\|^2\big).
\end{eqnarray}
Hence from equations $(\ref{am.1.14})$, $(\ref{am.1.22})$ and $(\ref{am.1.16})$, it follows that we end the proof of this theorem.
\end{pf}

We can take it to account that under a variation of a surface $M$, the $\mathsf{Area}(\varphi_t(M))=\int d\Omega_t$. Indeed, the volume element of the variation $d\Omega_t=u_td\Omega$, where $u_t$ is the area of the parallelogram spanned by $d(\varphi_t(X))$ and $d(\varphi_t(Y))$ for $X, Y\in \chi(M)$. Hereof, the variational vector field is decompose to $\xi=\xi^t+fN$, it is known that $\delta d\Omega_t= (\mathsf{div} \xi^t -2fH) d\Omega$, where $H$ denotes the mean curvature of $M$.

\begin{theorem}
Let $M$ be an isometrically immersed close surface in the $3-$dimensional harmonic conformally flat space $(\widetilde{M},\widetilde{g})$. Then $M$ is a minimal surface if the sign of mean curvature does not change.
\end{theorem}
\begin{pf}
Let $\varphi_t : M \rightarrow (\widetilde{M}^3,\widetilde{g})$ be a normal variation, where $\xi=fN$, of a close surface $M$ in the 3-dimensional harmonic conformally flat space $\widetilde{M}$. Then directly by applying the Gauss-Bonnet theory together with Theorem \ref{am.1.18} and Proposition \ref{am.1.10} we obtain
\begin{eqnarray}\label{am.1.26}
0=\delta\int_M K d\Omega &=& \int_M \delta K d\Omega + \int_M K(\delta d\Omega)\nonumber\\
&=& \int_M \big\{ \frac{f}{2}\big( (\lambda_1-\lambda_2) \widetilde{g}(\widetilde{\nabla}_{e_1}\omega^{\sharp}, e_1)+ (\lambda_2-\lambda_1)\widetilde{g}(\widetilde{\nabla}_{e_2} \omega^{\sharp},e_2)\big)\nonumber\\
&&+ \lambda_1\mathsf{Hess}_f(e_2,e_2)+\lambda_2\mathsf{Hess}_f(e_1,e_1)-\frac{1}{2}\mathsf{div}(\delta \omega^{\sharp})\nonumber\\
 &&+ f\big(2H(\lambda_1\lambda_2) - H\omega(\frac{\nabla f}{f}) -2HK +\frac{1}{8}H\| \omega^{\sharp}\|^2 - \omega(\frac{A\nabla f}{f})\big) \big\}d\Omega\nonumber\\
 &=& \int_M f\big \{ \lambda_1 \widetilde{g}(\widetilde{\nabla}_{e_1}\omega^{\sharp}, e_1)+ \lambda_2\widetilde{g}(\widetilde{\nabla}_{e_2} \omega^{\sharp},e_2)-\frac{3}{8}H\| \omega^{\sharp}\|^2- \omega(\frac{A\nabla f}{2f}) - H\omega(\frac{\nabla f}{f})\big\}d\Omega,
\end{eqnarray}
on the one hand we have
\begin{eqnarray}\label{am.1.32}
\lambda_2\widetilde{g}(\nabla_{e_2}\omega^{\sharp}, e_2)+\lambda_1\widetilde{g}(\nabla_{e_1}\omega^{\sharp},e_1) = 2H\mathsf{div}\omega^{\sharp} -\big(\lambda_2 \widetilde{g}(\widetilde{\nabla}_{e_1}\omega^{\sharp}, e_1)+ \lambda_1\widetilde{g}(\widetilde{\nabla}_{e_2} \omega^{\sharp},e_2)\big),
\end{eqnarray}
so we replace (\ref{am.1.32}) in $(\ref{am.1.26})$ which together with the Proposition \ref{am.1.23} follow
\begin{eqnarray}
  &&\int_M \big\{ -\big(\sigma\lambda_2 \widetilde{g}(\widetilde{\nabla}_{e_1}\nabla f, e_1)+ \sigma\lambda_1\widetilde{g}(\widetilde{\nabla}_{e_2} \nabla f,e_2)\big)+\frac{\sigma}{2}\omega(A\nabla f)\nonumber\\
 && + f\big(2H\mathsf{div}\omega^{\sharp} -\frac{7}{8}H\| \omega^{\sharp}\|^2-\omega(\frac{A\nabla f}{2f}) - H\omega(\frac{\nabla f}{f})\big) \big\}d\Omega=0.\label{am.1.27}
\end{eqnarray}
Hereof, $\widetilde{\nabla}_N N=-\frac{\nabla f}{f}$ \citet{V} and notice the assumption that 
\begin{eqnarray}\label{am.1.20}
0=\mathsf{Div}\omega^{\sharp}&=& \mathsf{div}\omega^{\sharp} +\widetilde{g}(\widetilde{\nabla}_{N}\omega^{\sharp}, N)\nonumber\\
&=&\mathsf{div}\omega^{\sharp}+\widetilde{g}(\omega^{\sharp},\frac{\nabla f}{f}),
\end{eqnarray}
it follows $\mathsf{div}\omega^{\sharp}=-\omega(\frac{\nabla f}{f})$ . From here, $(\ref{am.1.27})$ yields
\begin{eqnarray}\label{am.1.35}
&&\int_M \big\{ -f\big(\frac{7}{8}H\| \omega^{\sharp}\|^2+\omega(\frac{A\nabla f}{2f}) + 3H\omega(\frac{\nabla f}{f})\big)\nonumber\\
&&-\sigma\big(\lambda_2 \widetilde{g}(\widetilde{\nabla}_{e_1}\nabla f, e_1)+ \lambda_1\widetilde{g}(\widetilde{\nabla}_{e_2} \nabla f,e_2)\big) +\frac{\sigma}{2}\omega(A\nabla f)\big \}d\Omega=0.
\end{eqnarray}
After all, under the normal variation $\xi=N$, where $f\equiv 1$ the equation (\ref{am.1.35}) follows
\begin{eqnarray}
\int_M H \| \omega^{\sharp}\|^2 d\Omega=0,\nonumber
\end{eqnarray}
hence as respect to the assumption that the sign of mean curvature does not change, also $\omega^{\sharp}\neq 0$, it obtains that $H=0$ and we conclude the desired result.
\end{pf}
By regarding the mean curvature and considering the critical points of its integrating we get the following results.
\begin{theorem}
Let $\varphi:M \rightarrow \widetilde{M}^3$ be an isometric immersion of a close complete oriented surface $M$ in $3-$dimensional harmonic conformally flat space. If $M$ be a critical point of the mean curvature functional, that is, $\delta \int H d\Omega=0$, then $M$ is homeomorphic to a sphere and the Euler- lagrange equation associated to this functional satisfies
\begin{eqnarray}
  \| \omega^{\sharp}\|^2 -\lambda_1\lambda_2-K=0.
\end{eqnarray}
\end{theorem}
\begin{pf}
Let $\varphi_t: M \rightarrow \widetilde{M}^3$ denotes a normal variation of a close complete oriented surface $M$ in the 3-dimensional harmonic conformally flat space $\widetilde{M}^3$ with the variational vector field $\xi=fN$. Let $e_i$ for $i=1,2$ and $N$ be the principal directions corresponding to the eigenvalues $\lambda_i$ and local unit normal vector field on $M$, respectively. Then, from (\ref{am.1.2}) and Proposition (\ref{am.1.10}) by integrating of the mean curvature $H=\frac{\lambda_1+\lambda_2}{2}$ on $M$ we get
 \begin{eqnarray}\label{am.1.25}
\delta \int_M H d\Omega &=& \int_M \delta H d\Omega+ \int_M H (\delta d\Omega)\nonumber\\
&=& \frac{1}{2}\int_M \big\{ f\big(\widetilde{g}(\widetilde{R}(e_1, N)N, e_1)+ \widetilde{g}(\widetilde{R}(e_2, N)N, e_2) + \lambda_1^2 +\lambda_2^2 \big)\nonumber\\
&&+\mathsf{Hess}_f(e_1,e_1) + \mathsf{Hess}_f(e_2,e_2)-2fH^2 \big\} d\Omega,\nonumber\\
&=& \int_M f\big\{ -\lambda_1\lambda_2+ \frac{1}{2}\mathsf{\widetilde{Ric}}(N, N)\big\} d\Omega.
\end{eqnarray}
 where from $(\ref{am.1.1})$, the second terms of (\ref{am.1.25}) tensor $\mathsf{\widetilde{Ric}}$ associated to $\widetilde{M}$ satisfies
\begin{eqnarray}\label{am.1.28}
\mathsf{\widetilde{Ric}}(N, N) &=& -B(N, N)-\frac{1}{2}\big((\widetilde{g}(\widetilde{\nabla}_{e_1}\omega^{\sharp},
e_1)+\widetilde{g}(\widetilde{\nabla}_{e_2}\omega^{\sharp}, e_2)\big)+\frac{\| \omega^{\sharp}\|^2}{4},\nonumber\\
&=&-\frac{1}{2}\omega(\frac{\nabla f}{f})+\frac{3}{4}\|\omega^{\sharp}\|^2,
\end{eqnarray}
because it follows from $(\ref{am.1.21})$ and $(\ref{am.1.20})$ that $B(N,N)= \omega(\frac{\nabla f}{f})-\frac{\|\omega ^{\sharp}\|^2}{2}$ and  $\sum_{i=1}^2 \widetilde{g}(\widetilde{\nabla}_{e_i}\omega^{\sharp}, e_i)=\mathsf{div}\omega^{\sharp}= -\omega(\frac{\nabla f}{f})$, respectively. From here and (\ref{am.1.8}) we get
\begin{eqnarray}\label{am.1.34}
K= \frac{1}{2}\omega(\frac{\nabla f}{f})+\frac{1}{4}\|\omega^{\sharp}\|^2+\lambda_1\lambda_2,
\end{eqnarray}
 thus $(\ref{am.1.25})$ together with $(\ref{am.1.28})$ and $(\ref{am.1.34})$ imply
 \begin{eqnarray}\label{am.1.29}
 \delta \int_M H d\Omega = -\frac{1}{2}\int_M f\big( \lambda_1\lambda_2 +K - \|\omega^{\sharp}\|^2\big)d\Omega.
 \end{eqnarray}
According to the assumption $\delta \int_M H d\Omega=0$, hence $(\ref{am.1.29})$ follows that the critical points of the functional satisfy
\begin{eqnarray}
\lambda_1\lambda_2 +K - \|\omega^{\sharp}\|^2=0.
\end{eqnarray}
Further, equation (\ref{am.1.34}) yields $\lambda_1\lambda_2 = K - \frac{1}{2}\omega(\frac{\nabla f}{f})-\frac{1}{4}\|\omega^{\sharp}\|^2$ that together with $(\ref{am.1.29})$ imply
\begin{eqnarray}\label{am.1.36}
\delta \int_M H d\Omega= -\frac{1}{2}\int_M f\big(2K-\frac{1}{2}\omega(\frac{\nabla f}{f})-\frac{5}{4}\|\omega^{\sharp}\|^2\big)d\Omega,
\end{eqnarray}
such that under the normal variation $\xi=N$, that is, $f\equiv 1$ by taking into account the assumption and making use of the Gauss-Bonnet theory, equation $(\ref{am.1.36})$ follows
 \begin{eqnarray*}
  \chi(M)=\frac{5}{16\pi}\int \|\omega^{\sharp}\|^2 d\Omega >0,
 \end{eqnarray*}
 and we conclude the proof.
\end{pf}
Finding the critical points of a functional is a major problem in the calculus of variation. In the rest of this section for given a topological space we consider close surfaces in 3-dimensional harmonic conformally flat spaces, in which $\delta \int H^2 d\Omega=0$, that is known as the Willmore surfaces.
\begin{theorem}
Let $\varphi : M \rightarrow \widetilde{M}^3$ be an isometric immersion of a close surface $M$ in the 3-dimensional harmonic conformal flat space $\widetilde{M}^3$. Then $M$ is a Willmore surface if and only if
\begin{eqnarray}
\Delta H+ H\big( 2H^2- K + \| \omega^{\sharp}\|^2-\lambda_1\lambda_2\big)=0,
\end{eqnarray}
\end{theorem}
\begin{pf}
Analogously, we are in the hypotheses of the last theorem such that under a normal variation where $\xi=fN$ and equation $(\ref{am.1.2})$ we obtain
\begin{eqnarray}\label{am.1.33}
\delta \int_M H^2 d\Omega &=& \int_M 2H(\delta H) d\Omega+ \int_M H^2 (\delta d\Omega),\nonumber\\
&=& \int_M 2H\big (f\widetilde{K}+\frac{f}{2}\mathsf{\widetilde{Ric}}(N, N)+\frac{\Delta f}{2}+ f(2H^2-K)\big) -2fH^3 d\Omega,\nonumber\\
&=&\int_M fH \big(\|\omega^{\sharp}\|^2-\lambda_1\lambda_2-K \big)+H\Delta f+ 2fH^3 d\Omega,
\end{eqnarray}
where $\mathsf{\widetilde{Ric}(N,N)}$ satisfies $(\ref{am.1.28})$. Now, the divergence theorem can be applied such that equation (\ref{am.1.33}) becomes
\begin{eqnarray}
\delta \int H^2 d\Omega &=& \int f \big(\|\omega^{\sharp}\|^2H-\lambda_1\lambda_2 H +\Delta H+ H(2H^2-K)\big) d\Omega.
\end{eqnarray}
Consequently, the critical points of the above functional satisfies
\begin{eqnarray}
\Delta H+ H\big(\|\omega^{\sharp}\|^2-\lambda_1\lambda_2 +2H^2-K \big)=0,
\end{eqnarray}
thus we conclude the proof of result.
\end{pf}

\section{Acknowledge}
We will be really grateful the referees who help the paper to get much more improvement both in the english and mathematic.
\bibliographystyle{cas-model2-names}

\bibliography{cas-refs}


\end{document}